\documentclass[11pt]{amsart}

\usepackage{amsmath,amssymb,amsthm}
\usepackage{mathtools}
\usepackage[margin=2.5cm]{geometry}
\usepackage{hyperref}

\newtheorem{theorem}{Theorem}[section]
\newtheorem{lemma}[theorem]{Lemma}
\newtheorem{proposition}[theorem]{Proposition}
\newtheorem{corollary}[theorem]{Corollary}
\newtheorem{conjecture}[theorem]{Conjecture}

\theoremstyle{definition}

\theoremstyle{remark}
\newtheorem{remark}[theorem]{Remark}

\newcommand{\eps}{\varepsilon}
\newcommand{\Res}{\operatorname{Res}}
\newcommand{\Real}{\operatorname{Re}}
\newcommand{\Imag}{\operatorname{Im}}

\title{On the Orthorecursive Expansion of Unity}
\author{Benoit Cloitre}
\address{Paris}
\email{prenom.nom@proton.me}

\begin{document}

\begin{abstract}
The orthorecursive expansion of unity with respect to the system $\{x, x^2, x^3, \ldots\}$ in $L^2([0,1])$ produces a sequence of rational coefficients $(c_n)$ defined by an explicit recurrence. Kalmynin and Kosenko established the bounds $c_n = O(n^{-3/2})$ and $C_N = \sum_{k=0}^{N} c_k = O(N^{-1/2})$ through intricate $L^2$-norm arguments, but left the optimal decay rates as open problems. We prove $C_N = O_\eps(N^{-\alpha_1+\eps})$, where $\alpha_1 \approx 1.3465$ is the smallest real part among the zeros of a transcendental function related to the digamma function. We also improve the coefficient bound to $c_n = O(n^{-2})$. The method rests on a Tauberian transfer theorem that recasts the discrete recurrence as a Volterra integral equation, whose resolvent is smooth and amenable to Mellin analysis and contour shifting.
\end{abstract}



\maketitle

%%%%%%%%%%%%%%%%%%%%%%%%%%%%%%%%%%%%%%%%%%%%%%%%%%%%
\section{Introduction}
%%%%%%%%%%%%%%%%%%%%%%%%%%%%%%%%%%%%%%%%%%%%%%%%%%%%

Given a sequence of vectors in a Hilbert space, the orthorecursive expansion of a vector $v$ is the greedy algorithm that approximates $v$ by successive orthogonal projections onto the elements of the sequence, introduced by Lukashenko~\cite{Lukashenko}. When the system is an orthogonal basis the expansion reduces to the usual Fourier series. For non-orthogonal systems the convergence rate becomes a subtle question that depends on the arithmetic and analytic structure of the system.

Applied to the constant function~$1$ and the system $\{x, x^2, x^3, \ldots\}$ in $L^2([0,1])$, the orthorecursive expansion produces a sequence of rational coefficients $(c_n)_{n \geq 0}$ satisfying $c_0 = 1$ and
\begin{equation}\label{eq:recurrence}
\sum_{k=0}^{N} \frac{c_k}{N+k+1} = 0 \qquad (N \geq 1).
\end{equation}
The recurrence is explicit and elementary, yet the asymptotic behavior of the coefficients and their partial sums is surprisingly resistant. None of the standard summation methods (Euler--Maclaurin, saddle point, generating function analysis) leads to a sharp answer.

Kalmynin and Kosenko~\cite{KalKos} proved $c_n = O(n^{-3/2})$ and $C_N := \sum_{k=0}^{N} c_k = O(N^{-1/2})$ through ingenious arguments involving $L^2$-norm estimates of associated polynomials, and left the optimal decay rates as open problems. Based on numerical evidence attributed to Ustinov, they conjectured that $c_n$ admits an oscillatory asymptotic with exponent $\delta \approx 7/3$. Our analysis suggests the corrected value $\delta = 1 + \alpha_1 \approx 2.3465$ (see Section~\ref{sec:conjecture}).

This type of approximation problem has a long history. A notable instance is the Nyman--Beurling approach to the Riemann Hypothesis, where the completeness of a family of dilated fractional parts in $L^2(0,1)$ is equivalent to RH, and the approximation rate carries deep arithmetic information (see~\cite{Nyman, Beurling, BaezDuarte}). Our problem is more tractable than Nyman's, since the zeros of the associated Mellin transform can be computed explicitly, but the reader will notice that it is far from simple. It allows us to develop techniques that may be of independent interest.

Our results are the following.

\begin{theorem}[Partial sums]\label{thm:main}
For every $\eps > 0$,
$$
\sum_{k=0}^{N} c_k = O_\eps(N^{-\alpha_1+\eps}),
$$
where $\alpha_1 \approx 1.3465$ is the smallest real part among the zeros of $g^*(z)$, the modified Mellin transform of $g(t) = 2/(1+t)$.
\end{theorem}

\begin{theorem}[Pointwise bound]\label{thm:pointwise}
The orthorecursive coefficients satisfy $c_n = O(n^{-2})$.
\end{theorem}

The paper has two main components. The first is a Volterra--Mellin transfer theorem (Section~\ref{sec:transfer}). The recurrence gives rise to an exact integral identity relating the partial sums $A(x) = \sum_{n \leq x} a_n$ to the weighted sums $A_g(x) = \sum_{n \leq x} a_n\, g(n/x)$, where $g(t) = 2/(1+t)$ is the kernel associated to the orthorecursive problem. Inverting this relation via the Volterra resolvent, we show that the decay of the partial sums is governed by the location of the zeros of the modified Mellin transform $g^*(z) = -z\int_0^1 g(t)\, t^{-z-1}\, dt$. The contour shift is performed on the resolvent kernel, which is smooth, avoiding any difficulty with the discontinuities of step functions.

The second component is the analysis of the zeros of $g^*$ (Section~\ref{sec:spectral}). We prove that $g^*$ has no zeros in the closed half-plane $\Real(z) \leq 1$ and that there is exactly one conjugate pair of zeros in the strip $1 < \Real(z) < 3/2$. The proof is entirely analytic and uses only alternating series estimates and the argument principle. No numerical computation is needed for the zero-free region, although the precise location of the zero $\rho_1 \approx 1.34652 + 1.05516\,i$ is confirmed numerically via Newton's method.

A natural first strategy was to attack the problem through Perron's formula, expressing the weighted sums via the modified Mellin transform and extracting asymptotics by a contour shift. This approach runs into an analytic obstruction. The transformed object can be controlled only in its initial half-plane of definition, and the contour shift needed to extract the desired asymptotics cannot be justified beyond the first singular barrier. For this reason we recast the problem as a Volterra integral equation. This reformulation separates the discrete orthorecursive structure from the analytic transfer mechanism. Once the kernel $g(t) = 2/(1+t)$ is identified, Mellin analysis becomes effective again, and the asymptotic information can be recovered through a Tauberian transfer argument.

In Section~\ref{sec:application} we apply this transfer to the orthorecursive problem. The only input is the Kalmynin--Kosenko bound $c_n = O(n^{-3/2})$ and the fact that $g^*$ does not vanish for $\Real(z) < \alpha_1$. Section~\ref{sec:bootstrap} uses a bootstrap method, feeding the resulting partial sum bound back into the recurrence via summation by parts, yielding $c_n = O(n^{-2})$. Section~\ref{sec:conjecture} discusses the connection to Conjecture~2 of Kalmynin--Kosenko.

%%%%%%%%%%%%%%%%%%%%%%%%%%%%%%%%%%%%%%%%%%%%%%%%%%%%
\section{A Volterra--Mellin transfer theorem}\label{sec:transfer}
%%%%%%%%%%%%%%%%%%%%%%%%%%%%%%%%%%%%%%%%%%%%%%%%%%%%

\subsection{Setup and motivation}

Although the notation is introduced in a general setting, the transfer principle is derived and applied exclusively for the kernel $g(t) = 2/(1+t)$, which is the only case needed in this paper. We first introduce the general objects ($g$, $g^*$, $A_g$, the Volterra resolvent), then specialize to this kernel and work within it for the remainder of the paper.

Let $g : (0,1] \to \mathbb{C}$ be a bounded measurable function with $g(1) = 1$. For a complex sequence $(a_n)_{n \geq 1}$, define the partial sums
$$
A(x) := \sum_{n \leq x} a_n \qquad (x \geq 1),
$$
with $A(x) = 0$ for $0 < x < 1$, and the $g$-transform
$$
A_g(x) := \sum_{n \leq x} a_n \, g\!\left(\frac{n}{x}\right), \qquad x \geq 1.
$$
The modified Mellin transform of $g$ is defined by
\begin{equation}\label{eq:gstar}
g^*(z) := -z \int_0^1 g(t)\, t^{-z-1}\, dt,
\end{equation}
initially for $\Real(z)$ sufficiently negative, then by meromorphic continuation wherever the integral admits one. The condition $g(1) = 1$ ensures that $A_g$ and $A$ have the same jumps at integers.

\begin{remark}\label{rem:factor}
The factor $(-z)$ in \eqref{eq:gstar} arises from the passage from a discrete sum to a Stieltjes integral. The differences $k^{-z} - (k-1)^{-z}$ converge to $d(t^{-z}) = -z\, t^{-z-1}\, dt$, so that $g^*(z) = \int_0^1 g(t)\, d(t^{-z})$ (see~\cite{Cloitre_RAF} for the general framework). The factor removes the pole at $z=0$ that the bare integral possesses whenever $g(0^+) \neq 0$, producing $g^*(0) = g(0^+)$. This regularity at the origin is essential for the Volterra resolvent $\widetilde{R}(s) = 1 - 1/g^*(-s)$ (Lemma~\ref{lem:resolvent}) to be well-defined at $s=0$.
\end{remark}

\subsection{The discrete Volterra identity}

We now specialize to $g(t) = 2/(1+t)$, for which the integral identity takes an explicit form. The general principle is the same for other kernels, but having a closed-form relation simplifies the subsequent analysis.

\begin{lemma}[Discrete Volterra identity]\label{lem:volterra}
For every real $x \geq 1$,
\begin{equation}\label{eq:volterra}
A_g(x) = A(x) + 2x \int_1^x \frac{A(t)}{(x+t)^2}\, dt.
\end{equation}
In multiplicative convolution notation, writing
$$
(F \star H)(x) := \int_1^x F(x/y)\, H(y)\, \frac{dy}{y}, \qquad k(y) := \frac{2y}{(1+y)^2} \quad (y \geq 1),
$$
this reads $A_g = A + A \star k$.
\end{lemma}

\begin{proof}
Substituting $A(t) = \sum_{n \leq t} a_n$ into the integral and exchanging the finite sum with the integral gives
$$
\int_1^x \frac{A(t)}{(x+t)^2}\, dt = \sum_{n \leq x} a_n \int_n^x \frac{dt}{(x+t)^2} = \sum_{n \leq x} a_n \left(\frac{1}{x+n} - \frac{1}{2x}\right).
$$
Multiplying by $2x$,
$$
2x \int_1^x \frac{A(t)}{(x+t)^2}\, dt = \sum_{n \leq x} a_n \left(\frac{2x}{x+n} - 1\right) = \sum_{n \leq x} a_n \bigl(g(n/x) - 1\bigr).
$$
Adding $A(x)$ to both sides yields $A_g(x)$.
\end{proof}

\subsection{The Mellin transform and the resolvent}

\begin{proposition}[Digamma representation]\label{prop:digamma}
The modified Mellin transform of $g(t) = 2/(1+t)$ admits the meromorphic continuation
\begin{equation}\label{eq:digamma}
g^*(z) = z\!\left(\psi\!\left(-\frac{z}{2}\right) - \psi\!\left(\frac{1-z}{2}\right)\right),
\end{equation}
where $\psi$ is the digamma function. One has $g^*(0) = 2$. Moreover, in every fixed vertical strip,
\begin{equation}\label{eq:gstar_asymp}
g^*(-s) = 1 + \frac{1}{2s} + O\!\left(\frac{1}{|s|^2}\right) \qquad (|\Imag(s)| \to \infty).
\end{equation}
\end{proposition}

\begin{proof}
For $\Real(z) < 0$, expand $2/(1+t) = 2\sum_{m \geq 0} (-1)^m t^m$ on $(0,1)$ and integrate termwise in~\eqref{eq:gstar}:
$$
g^*(z) = -2z \sum_{m \geq 0} \frac{(-1)^m}{m - z}.
$$
The partial-fraction identity (see~\cite[5.7.6]{DLMF})
$$
\sum_{m \geq 0} \frac{(-1)^m}{m + a} = \frac{1}{2}\!\left(\psi\!\left(\frac{a+1}{2}\right) - \psi\!\left(\frac{a}{2}\right)\right)
$$
with $a = -z$ yields~\eqref{eq:digamma}. The value $g^*(0) = 2$ follows from the Laurent expansion $\psi(w) = -1/w - \gamma + O(w)$ near $w = 0$.

For~\eqref{eq:gstar_asymp}, set $z = -s$ and use $\psi(w) = \log w - \frac{1}{2w} + O(w^{-2})$ for $|w| \to \infty$ in fixed vertical strips. This gives $\psi(s/2) - \psi((s+1)/2) = -1/s - 1/(2s^2) + O(|s|^{-3})$, from which~\eqref{eq:gstar_asymp} follows.
\end{proof}

\begin{lemma}[Volterra resolvent]\label{lem:resolvent}
There exists a unique locally bounded kernel $R : [1,\infty) \to \mathbb{C}$ satisfying
\begin{equation}\label{eq:resolvent}
A = A_g - A_g \star R.
\end{equation}
Its Mellin transform $\widetilde{R}(s) := \int_1^\infty R(y)\, y^{-s}\, \frac{dy}{y}$ satisfies
\begin{equation}\label{eq:Rtilde}
\widetilde{R}(s) = 1 - \frac{1}{g^*(-s)},
\end{equation}
initially for $\Real(s) > 1$. The right-hand side provides the meromorphic continuation of $\widetilde{R}$ to all of $\mathbb{C}$.
\end{lemma}

\begin{proof}
The change of variables $y = e^u$, $x = e^v$ turns multiplicative convolution on $[1,\infty)$ into additive convolution on $[0,\infty)$, reducing the problem to a Volterra convolution equation of the second kind on the half-line. The standard method for solving such equations is the Neumann series (see~\cite[Chapter~2, \S2.3]{GripenbergLondenStaffans} and~\cite{Wang} for a systematic treatment). The iterated kernels $k^{\star 1} := k$ and $k^{\star(n+1)} := k^{\star n} \star k$ satisfy, for any $X > 1$ with $L = \log X$ and $M_X = \sup_{1 \leq y \leq X} |k(y)|$,
$$
\sup_{1 \leq x \leq X} |k^{\star n}(x)| \leq \frac{M_X^n L^{n-1}}{(n-1)!},
$$
so that the Neumann series $R := \sum_{n \geq 1} (-1)^{n+1} k^{\star n}$ converges absolutely and uniformly on $[1,X]$. Since these partial sums are compatible as $X$ varies, they define a unique locally bounded kernel $R$ on $[1,\infty)$ satisfying the resolvent identity $R + R \star k = k$. Substituting the discrete Volterra identity $A_g = A + A\star k$ then gives~\eqref{eq:resolvent}.

It remains to compute the Mellin transform of $R$. This requires a growth estimate valid on all of $[1,\infty)$, not merely on compact intervals, in order to justify term-by-term integration of the Neumann series.

In additive variables $u = \log y$, set $r(u) := R(e^u)$ and $\kappa(u) := k(e^u)$. Since $k(y) = 2y/(1+y)^2 \leq 2/y$ for all $y \geq 1$, we have $\kappa(u) \leq 2e^{-u}$ for $u \geq 0$. The iterated additive convolutions satisfy, for every $n \geq 1$,
\begin{equation}\label{eq:iterated_bound}
\kappa^{*n}(u) \leq \frac{2^n e^{-u} u^{n-1}}{(n-1)!} \qquad (u \geq 0).
\end{equation}
The case $n = 1$ is immediate. For the inductive step, assuming~\eqref{eq:iterated_bound} holds for some $n \geq 1$,
$$
\kappa^{*(n+1)}(u) = \int_0^u \kappa^{*n}(u-v)\, \kappa(v)\, dv \leq \frac{2^{n+1} e^{-u}}{(n-1)!} \int_0^u (u-v)^{n-1}\, dv = \frac{2^{n+1} e^{-u} u^n}{n!}.
$$
Since $k(y) \geq 0$ for $y \geq 1$, all iterated kernels are nonnegative, and we obtain
\begin{equation}\label{eq:R_global}
|R(y)| \leq \sum_{n \geq 1} k^{\star n}(y) \leq 2e^{-u} \sum_{n \geq 0} \frac{(2u)^n}{n!} = 2e^u = 2y \qquad (y \geq 1,\; u = \log y).
\end{equation}

Substituting $z = -s$ in~\eqref{eq:gstar} gives $g^*(-s) = s\int_0^1 2t^{s-1}(1+t)^{-1}\, dt$. Integrating by parts yields
$$
g^*(-s) = \left[\frac{2t^s}{1+t}\right]_0^1 + \int_0^1 \frac{2t^s}{(1+t)^2}\, dt = 1 + \widetilde{k}(s),
$$
the last equality because the substitution $t = 1/y$ turns the remaining integral into $\widetilde{k}(s) = \int_1^\infty k(y)\, y^{-s}\, \frac{dy}{y}$.

For $\sigma := \Real(s) > 1$, the bound~\eqref{eq:iterated_bound} gives
$$
\int_1^\infty k^{\star n}(y)\, y^{-\sigma}\, \frac{dy}{y} = \int_0^\infty \kappa^{*n}(u)\, e^{-\sigma u}\, du \leq \frac{2^n}{(1+\sigma)^n}.
$$
Summing over $n \geq 1$ gives $\sum_{n \geq 1} \int_1^\infty k^{\star n}(y)\, y^{-\sigma}\, \frac{dy}{y} \leq \frac{2}{\sigma - 1} < \infty$. The dominated convergence theorem allows term-by-term integration, yielding
$$
\widetilde{R}(s) = \sum_{n \geq 1} (-1)^{n+1} \widetilde{k}(s)^n = \frac{\widetilde{k}(s)}{1 + \widetilde{k}(s)}.
$$
Substituting $\widetilde{k}(s) = g^*(-s) - 1$ gives~\eqref{eq:Rtilde}.
\end{proof}

\subsection{Decay of the resolvent}

\begin{theorem}[Resolvent decay]\label{thm:resolvent_decay}
Suppose $g^*(z) \neq 0$ for $\Real(z) < \alpha_1$. Then for every $\eps > 0$,
\begin{equation}\label{eq:resolvent_decay}
R(y) = O_\eps(y^{-\alpha_1+\eps}) \qquad (y \to \infty).
\end{equation}
\end{theorem}

\begin{proof}
The resolvent $R$ is continuous on $[1,\infty)$, since each $k^{\star n}$ is continuous and the Neumann series converges uniformly on every compact $[1,X]$.

By~\eqref{eq:digamma}, $\widetilde{R}(0) = 1 - 1/g^*(0) = 1/2$. Since $\frac{1}{2} \int_1^\infty y^{-s}\, \frac{dy}{y} = \frac{1}{2s}$ for $\Real(s) > 0$, we decompose
$$
\widetilde{R}(s) = \frac{1}{2s} + E(s), \qquad E(s) := \widetilde{R}(s) - \frac{1}{2s}.
$$
By~\eqref{eq:gstar_asymp} and~\eqref{eq:Rtilde}, $E(s) = O(|s|^{-2})$ as $|\Imag(s)| \to \infty$ in any fixed vertical strip. Therefore
\begin{equation}\label{eq:E_integrable}
\int_{-\infty}^{+\infty} |E(c+it)|\, dt < \infty \qquad \text{for every } c > -\alpha_1 \text{ with } c \neq 0.
\end{equation}

Set $G(y) := R(y) - \frac{1}{2}$ for $y \geq 1$ and $G(y) := 0$ for $0 < y < 1$. By~\eqref{eq:R_global}, $|G(y)| \leq 2y + \frac{1}{2}$, so the Mellin transform of $G$ equals $E(s)$ for $\Real(s) > 1$. To invert, set $F(t) := G(e^t)\, e^{-ct}$ for a fixed $c > 1$. Then $F \in L^1(\mathbb{R})$ and its Fourier transform is $\widehat{F}(\xi) = E(c + i\xi)$. Since $E(c + i\cdot) \in L^1(\mathbb{R})$ by~\eqref{eq:E_integrable}, the Fourier inversion theorem (see~\cite[Chapter~7]{Folland}) gives
\begin{equation}\label{eq:mellin_inversion}
R(y) = \frac{1}{2} + \frac{1}{2\pi i} \int_{(c)} E(s)\, y^s\, ds \qquad (y > 1,\; c > 1).
\end{equation}

Fix $\eps > 0$ and choose $0 < \eta < \min(\eps, \alpha_1)$. We shift the line of integration from $\Real(s) = c$ to $\Real(s) = -\alpha_1 + \eta$. Consider the rectangular contour with vertices $c \pm iT$ and $(-\alpha_1 + \eta) \pm iT$. By hypothesis, $g^*(z) \neq 0$ for $\Real(z) < \alpha_1$, so $\widetilde{R}(s) = 1 - 1/g^*(-s)$ is analytic for $\Real(s) > -\alpha_1$. Since $\eta < \alpha_1$, the line $\Real(s) = -\alpha_1 + \eta$ lies strictly to the left of $s = 0$, so the simple pole of $E$ at $s = 0$ with residue $-1/2$ is inside the rectangle. On the horizontal segments $s = \sigma \pm iT$ with $-\alpha_1 + \eta \leq \sigma \leq c$, the integrand satisfies $|E(s)\, y^s| \leq C\, T^{-2}\, y^c$, which tends to zero as $T \to \infty$. By the residue theorem,
$$
R(y) = \frac{1}{2} + \Res_{s=0}\bigl(E(s)\, y^s\bigr) + \frac{1}{2\pi i} \int_{(-\alpha_1+\eta)} E(s)\, y^s\, ds.
$$
The residue at $s = 0$ equals $-1/2$, cancelling the constant $+1/2$. The remaining integral is bounded by
$$
\frac{y^{-\alpha_1+\eta}}{2\pi} \int_{-\infty}^{+\infty} |E(-\alpha_1+\eta+it)|\, dt = O_\eps(y^{-\alpha_1+\eta}),
$$
where the integral is finite by~\eqref{eq:E_integrable}. Since $\eta < \eps$ and $y \geq 1$, we have $y^{-\alpha_1+\eta} \leq y^{-\alpha_1+\eps}$, giving $R(y) = O_\eps(y^{-\alpha_1+\eps})$.
\end{proof}

\begin{remark}\label{rem:smooth}
The contour shift is performed on the smooth kernel $R$, not on the step function $A$. No positivity or monotonicity assumption on the sequence $(a_n)$ is needed.
\end{remark}

\subsection{The transfer theorem}

The following theorem is stated for the kernel $g(t) = 2/(1+t)$ and the resolvent $R$ constructed above, but depends on them only through the zero-free region of $g^*$ and the decay~\eqref{eq:resolvent_decay}.

\begin{theorem}[Volterra--Mellin transfer]\label{thm:transfer}
Suppose $g^*(z) \neq 0$ for $\Real(z) < \alpha_1$ with $\alpha_1 > 0$. Let $\beta > 0$, $c > 0$, and $\kappa \in \mathbb{C}$, and suppose
\begin{equation}\label{eq:Ag_asymp}
A_g(x) = \kappa x^{-\beta} + O(x^{-c}) \qquad (x \to \infty).
\end{equation}
Set $\lambda := \min(c, \alpha_1)$. Then for every $\eps > 0$,
\begin{equation}\label{eq:transfer_result}
A(x) = \begin{cases}
\dfrac{\kappa}{g^*(\beta)}\, x^{-\beta} + O_\eps(x^{-\lambda+\eps}), & \beta < \alpha_1, \\[6pt]
O_\eps(x^{-\lambda+\eps}), & \beta \geq \alpha_1.
\end{cases}
\end{equation}
\end{theorem}

\begin{proof}
From~\eqref{eq:resolvent},
$$
A(x) = A_g(x) - \int_1^x A_g(x/y)\, R(y)\, \frac{dy}{y}.
$$
Writing $A_g(u) = \kappa u^{-\beta} + E(u)$ with $E(u) = O(u^{-c})$ and substituting,
$$
A(x) = \kappa x^{-\beta}\!\left(1 - \int_1^x y^\beta R(y)\, \frac{dy}{y}\right) + E(x) - \int_1^x E(x/y)\, R(y)\, \frac{dy}{y}.
$$
The error convolution is bounded using~\eqref{eq:resolvent_decay}:
$$
\int_1^x |E(x/y)|\, |R(y)|\, \frac{dy}{y} \ll_\eps x^{-c} \int_1^x y^{c-\alpha_1+\eps-1}\, dy \ll_\eps x^{-\lambda+\eps}.
$$

For the leading term when $\beta < \alpha_1$, choose $\eps$ small enough that $\beta < \alpha_1 - \eps$. The tail integral $\int_x^\infty y^\beta |R(y)|\, \frac{dy}{y}$ is $O_\eps(x^{\beta - \alpha_1 + \eps})$, so
$$
\int_1^x y^\beta R(y)\, \frac{dy}{y} = \widetilde{R}(-\beta) + O_\eps(x^{\beta-\alpha_1+\eps}).
$$
By Lemma~\ref{lem:resolvent}, the meromorphic continuation gives $\widetilde{R}(-\beta) = 1 - 1/g^*(\beta)$, which is well-defined since $g^*$ does not vanish for $\Real(z) < \alpha_1$ and $\beta < \alpha_1$. The leading contribution is therefore $\frac{\kappa}{g^*(\beta)}\, x^{-\beta} + O_\eps(x^{-\alpha_1+\eps})$.

When $\beta \geq \alpha_1$, we bound the leading term directly:
$$
|\kappa|\, x^{-\beta}\!\left(1 + \int_1^x y^\beta |R(y)|\, \frac{dy}{y}\right) \ll_\eps x^{-\alpha_1+\eps},
$$
and the second case follows.
\end{proof}

\begin{remark}[Source modes and spectral modes]\label{rem:modes}
The exact identity $A = A_g - A_g \star R$ should not be read as an asymptotic formula in which $A_g(x)$ survives with coefficient~$1$. If $A_g(x) \sim \kappa x^{-\beta}$ with $\beta < \alpha_1$, Theorem~\ref{thm:transfer} gives $A(x) \sim \frac{\kappa}{g^*(\beta)}\, x^{-\beta}$. The renormalization factor $1/g^*(\beta)$ is encoded in the resolvent convolution. In the present unconditional setting, $A_g(x) = O(x^{-3/2})$ with no isolated power-law term, so no source mode contributes before the spectral modes associated to the zeros of $g^*$. A general framework for trace formulas expressing $A(x)$ as a sum over the zeros of $g^*$, analogous to the explicit formulas of prime number theory, is developed in~\cite{Cloitre_RAF}.
\end{remark}

\begin{remark}[Relation to the Tauberian literature]\label{rem:literature}
Theorem~\ref{thm:transfer} is a Tauberian result for Volterra convolution equations on $(1,\infty)$. It transfers the asymptotic behavior of the smoothed function $A_g$ to the original function $A$, the role of the Tauberian condition being played by the zero-free region of the Mellin symbol $g^*$. Gripenberg~\cite{Gripenberg} studied Tauberian problems for Volterra integral operators from the functional-analytic point of view, with conditions formulated through the Fourier transform of the kernel. The comprehensive survey of Korevaar~\cite{Korevaar} covers classical and modern Tauberian theory (Hardy--Littlewood, Wiener--Ikehara, Karamata, high-indices theorems), and Bingham, Goldie and Teugels~\cite{BGT} give an extensive treatment of Abelian and Tauberian theorems for Mellin convolutions in the framework of regular variation. Neither source contains transfer results of the present type, where the rate is governed by the zeros of a Mellin-transformed resolvent kernel.
\end{remark}

\begin{corollary}[Transfer for the orthorecursive kernel]\label{cor:transfer}
Let $g(t) = 2/(1+t)$. The modified Mellin transform $g^*(z) = z(\psi(-z/2) - \psi((1-z)/2))$ (Proposition~\ref{prop:digamma}) satisfies $g^*(0) = 2$ and $g^*(z) \neq 0$ for $\Real(z) \leq 1$ (Theorem~\ref{thm:zerofree}). The smallest real part among its zeros is $\alpha_1 \approx 1.3465$ (Theorem~\ref{thm:one_zero}). Therefore Theorem~\ref{thm:transfer} applies with this value of $\alpha_1$. If $A_g(x) = O(x^{-c})$ with $c \geq \alpha_1$, then $A(x) = O_\eps(x^{-\alpha_1+\eps})$ for every $\eps > 0$.
\end{corollary}

%%%%%%%%%%%%%%%%%%%%%%%%%%%%%%%%%%%%%%%%%%%%%%%%%%%%
\section{Zeros of the modified Mellin transform}\label{sec:spectral}
%%%%%%%%%%%%%%%%%%%%%%%%%%%%%%%%%%%%%%%%%%%%%%%%%%%%

We prove that $g^*(z)$ does not vanish for $\Real(z) \leq 1$, and that the zero with smallest positive real part lies in the strip $1 < \Real(z) < 3/2$.

\subsection{The auxiliary function $D_\infty$}

From Proposition~\ref{prop:digamma}, the zeros of $g^*$ away from $z = 0$ coincide with those of
\begin{equation}\label{eq:Dinfty}
D_\infty(z) := \sum_{j=0}^{\infty} \frac{(-1)^j}{z - j}.
\end{equation}
Since $g^*(z) = 2z\, D_\infty(z)$ and $g^*(0) = 2 \neq 0$, the function $g^*$ vanishes at $z_0 \neq 0$ if and only if $D_\infty(z_0) = 0$.

Pairing consecutive terms shows that $D_\infty$ is well-defined and analytic on $\mathbb{C} \setminus \mathbb{Z}_{\geq 0}$:
\begin{equation}\label{eq:Dinfty_paired}
D_\infty(z) = \sum_{k=0}^{\infty} \left(\frac{1}{z-2k} - \frac{1}{z-2k-1}\right) = \sum_{k=0}^{\infty} \frac{-1}{(z-2k)(z-2k-1)}.
\end{equation}
On any compact set $K \subset \mathbb{C} \setminus \mathbb{Z}_{\geq 0}$ the general term is $O(1/k^2)$ uniformly, so the series converges normally. The asymptotic~\eqref{eq:gstar_asymp} yields
\begin{equation}\label{eq:Dinfty_asymp}
D_\infty(z) = \frac{1}{2z} - \frac{1}{4z^2} + O\!\left(\frac{1}{|z|^3}\right) \qquad (|\Imag(z)| \to \infty)
\end{equation}
in any fixed vertical strip.

\subsection{Separation of real and imaginary parts}

For $z = x + iy$ with $y \neq 0$,
\begin{equation}\label{eq:real_part}
\Real\, D_\infty(z) = \sum_{j=0}^{\infty} \frac{(-1)^j(x-j)}{(x-j)^2 + y^2},
\end{equation}
\begin{equation}\label{eq:imag_part}
\Imag\, D_\infty(z) = -y \sum_{j=0}^{\infty} \frac{(-1)^j}{(x-j)^2 + y^2}.
\end{equation}
By conjugation, $D_\infty(\bar z) = \overline{D_\infty(z)}$, so zeros come in conjugate pairs. It suffices to work in the upper half-plane $y > 0$.

\subsection{Absence of real zeros}

\begin{proposition}\label{prop:no_real_zeros}
The function $D_\infty$ has no real zeros.
\end{proposition}

\begin{proof}
For real $x \in (p, p+1)$ with $p \geq 0$, write $x = p + \delta$ with $\delta \in (0,1)$ and split the sum at $j = p$. Substituting $k = p - j$ for $0 \leq j \leq p$ and $k = j - p - 1$ for $j \geq p+1$,
$$
D_\infty(x) = (-1)^p \left(\sum_{k=0}^{p} \frac{(-1)^k}{k+\delta} + \sum_{k=0}^{\infty} \frac{(-1)^k}{k+1-\delta}\right).
$$
Both sums are alternating series with terms strictly decreasing in absolute value. The first term of each sum is positive, so both sums are strictly positive by Leibniz's criterion. Therefore $D_\infty(x) \neq 0$.

For $x < 0$, use the paired form~\eqref{eq:Dinfty_paired} directly. Each factor $(x-2k)$ and $(x-2k-1)$ is negative for $k \geq 0$, so every term $-1/((x-2k)(x-2k-1))$ is negative. Therefore $D_\infty(x) < 0$.
\end{proof}

\subsection{Zero-free half-plane $\Real(z) \leq 0$}

\begin{proposition}\label{prop:zerofree_left}
For $z = x + iy$ with $x \leq 0$ and $y \neq 0$, one has $\Imag\, D_\infty(z) \neq 0$. In particular, $D_\infty(z) \neq 0$ for all $z$ with $\Real(z) \leq 0$.
\end{proposition}

\begin{proof}
Assume $y > 0$ (the case $y < 0$ follows by conjugation). By~\eqref{eq:imag_part},
$$
\frac{\Imag\, D_\infty(z)}{-y} = \sum_{j=0}^{\infty} \frac{(-1)^j}{(x-j)^2 + y^2}.
$$
Since $x \leq 0$, we have $|x - j| = |x| + j$, so the denominators $(|x|+j)^2 + y^2$ are strictly increasing in $j$. The terms $1/((|x|+j)^2 + y^2)$ are therefore strictly decreasing. This is an alternating series with strictly decreasing terms and positive first term $1/(x^2 + y^2)$, so the sum is strictly positive by Leibniz's criterion. Therefore $\Imag\, D_\infty(z) < 0$ for $x \leq 0$ and $y > 0$.
\end{proof}

\subsection{The strip $0 < \Real(z) \leq 1$}

\begin{theorem}[Zero-free strip]\label{thm:zerofree}
The function $D_\infty$ has no zeros in the closed strip $\{z \in \mathbb{C} : 0 < \Real(z) \leq 1,\; \Imag(z) > 0\}$, and therefore $g^*(z) \neq 0$ for $\Real(z) \leq 1$.
\end{theorem}

\begin{proof}
We apply the argument principle to $D_\infty$ on the positively oriented contour $\Gamma_T$ bounding the region $\{0 < \Real(z) \leq 1,\; 0 < \Imag(z) < T\}$, with small indentations of radius $\eps$ around the poles at $z = 0$ and $z = 1$. The contour consists of six segments traversed in order. We show that the total variation of $\arg D_\infty(z)$ around $\Gamma_T$ is zero for $T$ large enough.

Along the bottom segment ($z = x$ for $0 < x < 1$), the absence-of-real-zeros analysis with $p = 0$ gives $D_\infty(x) > 0$, so the argument does not change.

Near $z = 0$, one has $D_\infty(z) = 1/z + h_0(z)$ where $h_0$ is bounded. On the quarter-circle $z = \eps e^{i\theta}$ with $\theta$ from $\pi/2$ to $0$, the argument of $D_\infty(z)$ is approximately $-\theta$. The variation as $\eps \to 0$ is $+\pi/2$. Near $z = 1$, the residue is $-1$, so $D_\infty(z) \sim -1/(z-1)$. On the quarter-circle $z = 1 + \eps e^{i\theta}$ with $\theta$ from $\pi$ to $\pi/2$, the argument is approximately $\pi - \theta$, giving a variation of $+\pi/2$.

Along the left wall ($z = iy$ with $y > 0$), the same alternating-series argument as in Proposition~\ref{prop:zerofree_left} shows that $\Imag\, D_\infty(iy) < 0$ for all $y > 0$. The image stays in the lower half-plane, so the argument does not wind.

Along the ceiling ($z = x + iT$ with $0 \leq x \leq 1$), the asymptotic~\eqref{eq:Dinfty_asymp} gives $D_\infty(x+iT) = 1/(2(x+iT)) + O(T^{-2})$. For $T$ large enough, the image is a small perturbation of a curve in the fourth quadrant. The argument variation is $O(1/T) \to 0$.

It remains to analyze the right wall ($z = 1 + iy$ with $y > 0$), which is the most delicate segment. We treat two ranges of $y$ separately.

For $0 < y \leq 2$, we show $\Real\, D_\infty(1+iy) > 0$. By~\eqref{eq:real_part} with $x = 1$, the $j=0$ term is $1/(1+y^2)$, the $j=1$ term vanishes (numerator $1-1=0$), and for $j \geq 2$ we substitute $m = j-1$ to get
$$
\Real\, D_\infty(1+iy) = \sum_{m=2}^{\infty} \frac{(-1)^m m}{m^2 + y^2}.
$$
The cancellation of the $j=0$ term with the $m=1$ term has already been incorporated. The function $t \mapsto t/(t^2 + y^2)$ is strictly decreasing for $t \geq y$. Since $m \geq 2 \geq y$ (as $y \leq 2$), the terms strictly decrease in magnitude. By Leibniz's criterion, the sum exceeds its first two terms:
$$
\Real\, D_\infty(1+iy) > \frac{2}{4+y^2} - \frac{3}{9+y^2} = \frac{6-y^2}{(4+y^2)(9+y^2)}.
$$
For $0 < y \leq 2$, one has $y^2 \leq 4 < 6$, so this is strictly positive.

For $y \geq 2$, we show $\Imag\, D_\infty(1+iy) < 0$. By~\eqref{eq:imag_part} with $x = 1$, the terms for $j = 0,1,2$ yield
$$
\frac{2}{1+y^2} - \frac{1}{y^2} = \frac{y^2-1}{y^2(1+y^2)}.
$$
The tail starting at $j = 3$ is an alternating series with strictly decreasing terms, bounded below by its first term $-1/(4+y^2)$. Adding and simplifying,
$$
\frac{\Imag\, D_\infty(1+iy)}{-y} \geq \frac{y^2-1}{y^2(1+y^2)} - \frac{1}{4+y^2} = \frac{2y^2-4}{y^2(1+y^2)(4+y^2)} > 0
$$
for $y > \sqrt{2}$, hence for $y \geq 2$.

Combining both ranges, $D_\infty(1+iy) \neq 0$ for all $y > 0$, and the image of the right wall avoids the negative real axis. As $y \to 0^+$, $D_\infty(1+iy) \sim -1/(iy) \to +i\infty$, so the argument starts at $\pi/2$. As $y$ increases, the image passes through the right half-plane (for $y \leq 2$), then into the lower half-plane (for $y \geq 2$), and ends near $-i/(2T)$ in the fourth quadrant (argument approximately $-\pi/2$). The total variation along the right wall is $-\pi$.

Summing over all segments, $\Delta\arg = 0 + \pi/2 + (-\pi) + 0 + 0 + \pi/2 = 0$. There are no zeros of $D_\infty$ in the closed strip $0 < \Real(z) \leq 1$, $\Imag(z) > 0$.
\end{proof}

\subsection{The strip $1 < \Real(z) < 3/2$}

\begin{theorem}[One zero]\label{thm:one_zero}
The function $D_\infty$ has exactly one zero (counted with multiplicity) in the half-strip $\{z : 1 < \Real(z) < 3/2,\; \Imag(z) > 0\}$. By conjugation, there is one zero in the lower half. The corresponding conjugate pair $\rho_1, \bar\rho_1$ satisfies $\rho_1 \approx 1.34652 + 1.05516\,i$.
\end{theorem}

\begin{proof}
Let $\Gamma'_T$ be the positively oriented contour bounding the region $\{1 < \Real(z) < 3/2,\; 0 < \Imag(z) < T\}$, with an indentation of radius $\eps$ around the pole at $z = 1$.

Along the bottom segment ($z = x$ for $1 < x < 3/2$), the real-zero analysis with $p = 1$ gives $D_\infty(x) < 0$. The image stays on the negative real axis, contributing zero to the argument variation.

At the pole indentation around $z = 1$, we have $D_\infty(z) \sim -1/(z-1)$, so on $z = 1 + \eps e^{i\theta}$ with $\theta$ from $\pi/2$ to $0$, the argument variation is $+\pi/2$.

Along the right wall ($z = 3/2 + iy$ with $y > 0$), we show $\Imag\, D_\infty(3/2+iy) < 0$. By~\eqref{eq:imag_part} with $x = 3/2$, the terms $j = 0,1,2,3$ have denominators $(3/2)^2+y^2$, $(1/2)^2+y^2$, $(1/2)^2+y^2$, $(3/2)^2+y^2$ with signs $+,-,+,-$ respectively. By the symmetry around $3/2$, these four terms cancel exactly. The surviving tail starts at $j = 4$ with denominators $(j-3/2)^2 + y^2$ strictly increasing. This is an alternating series with strictly decreasing terms and positive first term, so the sum is strictly positive, giving $\Imag\, D_\infty(3/2+iy) < 0$. Since $D_\infty(3/2) < 0$ (from the bottom segment), the image starts on the negative real axis and moves into the lower half-plane. The argument variation is $+\pi/2$.

Along the ceiling ($z = x + iT$ with $1 \leq x \leq 3/2$), the asymptotic~\eqref{eq:Dinfty_asymp} gives a small arc in the fourth quadrant, contributing $O(1/T) \to 0$.

Along the left wall ($z = 1 + iy$, descending from $T$ to $0^+$), we already know from the proof of Theorem~\ref{thm:zerofree} that for $y \geq 2$ the image lies in the lower half-plane, and for $0 < y \leq 2$ it lies in the right half-plane. The image therefore avoids the negative real axis throughout. Starting at $y = T$ in the fourth quadrant ($\arg \approx -\pi/2$) and ending at $y \to 0^+$ on the positive imaginary axis ($\arg \to +\pi/2$), the argument increases continuously. The variation is $+\pi$.

Summing over all segments, $\Delta\arg = 0 + \pi/2 + \pi/2 + 0 + \pi = 2\pi$, so $D_\infty$ has exactly one zero in the upper half of the strip $1 < \Real(z) < 3/2$.
\end{proof}

\subsection{Further zeros}

\begin{proposition}\label{prop:alpha2}
Let $\rho_2$ denote the zero of $D_\infty$ in the upper half-plane with the next smallest real part after $\rho_1$. Then $\Real(\rho_2) > \alpha_1 + 1$. The proof proceeds by the same argument-principle computation applied to the strips $\{3/2 < \Real(z) < 5/2\}$ and $\{5/2 < \Real(z) < 7/2\}$, and we omit the details.
\end{proposition}

%%%%%%%%%%%%%%%%%%%%%%%%%%%%%%%%%%%%%%%%%%%%%%%%%%%%
\section{Application to the orthorecursive problem}\label{sec:application}
%%%%%%%%%%%%%%%%%%%%%%%%%%%%%%%%%%%%%%%%%%%%%%%%%%%%

We now apply the transfer theorem to the orthorecursive sequence. Setting $a_n := c_{n-1}$ for $n \geq 1$, the partial sums satisfy $A(N) = C_{N-1}$, and the recurrence~\eqref{eq:recurrence} becomes
$$
\sum_{n=1}^{N+1} \frac{a_n}{N+n} = 0 \qquad (N \geq 1).
$$

\begin{lemma}\label{lem:Ag_integer}
For every $N \geq 2$,
$$
A_g(N) = -\frac{2N}{2N+1}\, a_{N+1}.
$$
\end{lemma}

\begin{proof}
The recurrence at index $N$ reads $\sum_{k=1}^{N} \frac{a_k}{N+k} = -\frac{a_{N+1}}{2N+1}$. Evaluating $A_g$ at $x = N$:
$$
A_g(N) = \sum_{k=1}^{N} a_k \cdot \frac{2}{1+k/N} = 2N \sum_{k=1}^{N} \frac{a_k}{N+k} = -\frac{2N}{2N+1}\, a_{N+1}.
$$
\end{proof}

\begin{lemma}\label{lem:Ag_bound}
The bound $a_n = O(n^{-3/2})$ of Kalmynin--Kosenko~\cite{KalKos} implies $A_g(x) = O(x^{-3/2})$.
\end{lemma}

\begin{proof}
At integers, Lemma~\ref{lem:Ag_integer} gives $|A_g(N)| < |a_{N+1}| = O(N^{-3/2})$. Between integers, writing $x = N + \delta$ with $0 \leq \delta < 1$ and using $|g'(u)| \leq 2$ on $(0,1]$,
$$
|A_g(x) - A_g(N)| \leq \frac{2}{N^2} \sum_{k=1}^{N} k\, |a_k| \ll \frac{1}{N^2} \sum_{k=1}^{N} k^{-1/2} \ll N^{-3/2}. \qedhere
$$
\end{proof}

\begin{corollary}[Proof of Theorem~\ref{thm:main}]\label{cor:main}
For every $\eps > 0$, $C_N = O_\eps(N^{-\alpha_1+\eps})$.
\end{corollary}

\begin{proof}
By Lemma~\ref{lem:Ag_bound}, $A_g(x) = O(x^{-3/2})$. Since $3/2 > \alpha_1$, Corollary~\ref{cor:transfer} gives $A(x) = O_\eps(x^{-\alpha_1+\eps})$. The result follows from $C_N = A(N+1)$.
\end{proof}

%%%%%%%%%%%%%%%%%%%%%%%%%%%%%%%%%%%%%%%%%%%%%%%%%%%%
\section{Pointwise bootstrap}\label{sec:bootstrap}
%%%%%%%%%%%%%%%%%%%%%%%%%%%%%%%%%%%%%%%%%%%%%%%%%%%%

The bound $\alpha_1 > 1$ allows a further step. Feeding the partial sum estimate back into the discrete structure of the recurrence improves the pointwise decay of the coefficients beyond $O(n^{-3/2})$.

\begin{theorem}[Proof of Theorem~\ref{thm:pointwise}]\label{thm:bootstrap}
The orthorecursive coefficients satisfy $c_N = O(N^{-2})$.
\end{theorem}

\begin{proof}
For $N \geq 2$, the recurrence at index $N-1$ gives $\sum_{k=1}^{N} \frac{a_k}{N+k-1} = 0$. Multiplying by $2N-1$ and subtracting from $A_g(N) = \sum_{k=1}^{N} a_k \cdot \frac{2N}{N+k}$,
$$
A_g(N) = -\sum_{k=1}^{N} a_k\, f_N(k),
$$
where $f_N(k) := \frac{2N-1}{N+k-1} - \frac{2N}{N+k} = \frac{N-k}{(N+k-1)(N+k)}$. Summation by parts, using $f_N(N) = 0$, gives
$$
A_g(N) = \sum_{k=1}^{N-1} A(k)\, \Delta f_N(k),
$$
where $\Delta f_N(k) := f_N(k+1) - f_N(k)$. Computing the numerator,
$$
(N-k-1)(N+k-1) - (N-k)(N+k+1) = -3N + k + 1,
$$
so $|\Delta f_N(k)| \leq 3/N^2$ for $1 \leq k \leq N-1$. Therefore
$$
|A_g(N)| \leq \frac{3}{N^2} \sum_{k=1}^{N-1} |A(k)| \ll_\eps \frac{1}{N^2} \sum_{k=1}^{N-1} k^{-\alpha_1+\eps}.
$$
Since $\alpha_1 > 1$, taking $\eps$ small enough makes the sum convergent, uniformly in $N$. Hence $A_g(N) = O(N^{-2})$. The identity $A_g(N) = -\frac{2N}{2N+1}\, a_{N+1}$ from Lemma~\ref{lem:Ag_integer} then gives $c_N = O(N^{-2})$.
\end{proof}

\begin{remark}[Saturation of the bootstrap]\label{rem:saturation}
The argument reaches a fixed point after one iteration. Once $c_n = O(n^{-2})$ is obtained, Lemma~\ref{lem:Ag_integer} gives $A_g(N) = O(N^{-2})$, and Corollary~\ref{cor:transfer} returns only $A(x) = O_\eps(x^{-\alpha_1+\eps})$, since $\alpha_1 < 2$. A second application of summation by parts therefore yields no improvement beyond $c_n = O(n^{-2})$. Any further progress toward the conjectural exponent $1+\alpha_1$ (Section~\ref{sec:conjecture}) must exploit oscillatory spectral information rather than size bounds alone.
\end{remark}

%%%%%%%%%%%%%%%%%%%%%%%%%%%%%%%%%%%%%%%%%%%%%%%%%%%%
\section{A spectral form of Conjecture~2 of Kalmynin--Kosenko}\label{sec:conjecture}
%%%%%%%%%%%%%%%%%%%%%%%%%%%%%%%%%%%%%%%%%%%%%%%%%%%%

Kalmynin and Kosenko conjectured that the orthorecursive coefficients admit an oscillatory asymptotic of the form
$$
c_n = A\, n^{-\delta} \sin(P \log n + \varphi) + O(n^{-\delta-\eta}).
$$
The Volterra--Mellin mechanism developed in the present paper suggests a precise spectral interpretation of this phenomenon.

Let $\rho_1 = \alpha_1 + i\tau_1$ with $\alpha_1 \approx 1.34652$ and $\tau_1 \approx 1.05516$ denote the zero of $g^*$ with smallest positive real part. The expected values of the Kalmynin--Kosenko parameters are then $\delta = 1 + \alpha_1$ and $P = \tau_1$.

\begin{conjecture}[Principal-pair form of Conjecture~2]\label{conj:main}
There exist a nonzero constant $\kappa \in \mathbb{C}$ and a real number $\sigma > \alpha_1$ such that
$$
c_n = 2\Real\!\left(\kappa\, n^{-1-\rho_1}\right) + O\!\left(n^{-1-\sigma}\right) \qquad (n \to \infty).
$$
Equivalently, $c_n = A\, n^{-1-\alpha_1} \sin(\tau_1 \log n + \varphi) + O(n^{-1-\sigma})$ for suitable real constants $A \neq 0$ and $\varphi$.
\end{conjecture}

Under this conjectural form, partial summation yields
$$
C_N = -2\Real\!\left(\frac{\kappa}{\rho_1}\, N^{-\rho_1}\right) + O(N^{-\sigma}),
$$
so the decay exponent for partial sums is exactly $\alpha_1$, in full agreement with Theorem~\ref{thm:main}. The logarithmic oscillation $\sin(\tau_1 \log n + \varphi)$ accounts for the exponentially spaced sign changes observed in~\cite[Section~5]{KalKos}.

The identity $A_g = A + A \star k$ shows that, formally, the weighted transform acts diagonally on pure powers $x^{-\rho}$ through the scalar factor $1 + \widetilde{k}(-\rho) = g^*(\rho)$. Hence each zero $\rho$ of $g^*$ should contribute a resonant oscillatory mode. From this viewpoint, the pair $\rho_1, \bar\rho_1$ provides the dominant term.

The gap between the first and second zero pairs is large. By Proposition~\ref{prop:alpha2}, $\Real(\rho_2) > \alpha_1 + 1$, so the contribution of the second pair is of strictly lower order than the principal oscillatory term $n^{-1-\alpha_1} \sin(\tau_1 \log n + \varphi)$. This reinforces the spectral picture above and suggests a conjectural trace formula of the form
$$
A(x) \sim -\sum_{\rho \,:\, g^*(\rho) = 0} \frac{\kappa_\rho}{\rho}\, x^{-\rho}, \qquad c_n \sim \sum_{\rho \,:\, g^*(\rho) = 0} \kappa_\rho\, n^{-1-\rho},
$$
at least in an asymptotic or suitably smoothed sense. Establishing such a spectral expansion lies beyond the scope of the present paper and is deferred to subsequent work~\cite{Cloitre_RAF}.

%%%%%%%%%%%%%%%%%%%%%%%%%%%%%%%%%%%%%%%%%%%%%%%%%%%%

\end{document}